\documentclass[11pt]{article}
\usepackage{amssymb}
\usepackage{amsmath}
\usepackage{amscd}
\usepackage{latexsym}
\usepackage{graphicx}
\usepackage{url}
\usepackage{epsfig} 

\catcode`@=11
\makeatletter\renewcommand{\section}{\@startsection
{section}{1}{\z@}{-3.5ex plus -1ex minus
    -.2ex}{2.3ex plus .2ex}{\large\bf }}

\makeatletter\renewcommand{\subsection}{\@startsection{subsection}{2}{\z@}{-3.25ex
plus -1ex minus
   -.2ex}{1.5ex plus .2ex}{\bf }}

\numberwithin{equation}{section}
\newcounter{saveeqn}

\catcode`@=12
\def\bea{\begin{eqnarray}}
\def\eea{\end{eqnarray}}

\def\1{{\bar 1}}
\def\2{{\bar 2}}
\def\3{{\bar 3}}
\def\4{{\bar 4}}

\def\+{\dagger}
\def\={\ =\ }

\def\and{\quad\textrm{and}\quad}

\begin{document}

\begin{titlepage}
\setcounter{page}{0}

\vskip 2.0cm

\begin{center}

{\Large\bf
Computing the Mertens function on a GPU
}

\vspace{12mm}

{\large Eugene Kuznetsov}

\vspace{12mm}

\begin{abstract}
\noindent

A GPU implementation of an algorithm to compute the Mertens function in $O(x^{2/3+\epsilon})$ time is discussed. Results for x up to $10^{22}$, and a new extreme value for $M(x)/x^{1/2}$, -0.585768 ($M(x)\approx -1.996*10^9$ at $x \approx 1.161*10^{19}$), are reported.An approximate algorithm is used to examine values of $M(x)$ for $x$ up to $\exp{(10^{15})}$.

\end{abstract}

\end{center}
\end{titlepage}

\section{Introduction} 

For any positive integer x, ''Mobius function'' $\mu(x)$ is defined as:

* 0, if x is divisible by a square of any prime

* 1, if the factorization of x contains an even number of distinct primes

* -1, otherwise

Mertens function of a positive integer x is defined as the sum of values of the Mobius function for 1 to x:

$$M(n) = \sum_{x=1}^n \mu(x)$$

This function scales roughly as $n^{1/2}$ over the interval where it's been studied. An old hypothesis, "Mertens conjecture" (Stieltjes, 1885) proposed that $|M(n)| < n^{1/2}$ for all n. This was disproved by Odlyzko and te Riele (1985)~\cite{1}, though no explicit counterexample is known.

The largest known value of $|M(n)/n^{1/2}|$, 0.570591 at n=7766842813, was found in ~\cite{2}. Values of M(n) for all consecutive n up to $10^{14}$ were computed in ~\cite{3}. No larger values were found. 

A computer algorithm to calculate isolated values of M(n) in $O(n^{2/3+\epsilon})$ time was described in ~\cite{4} and used to calculate several values up to $10^{16}$.

The algorithm described next is, in a sense, a variation of the algorithm by ~\cite{4}. The running time is still $O(n^{2/3+\epsilon})$. A distinctive feature of this algorithm is that executing it at any value of n allows to simultaneously compute all values of $M(y)$ for integer $y=\lfloor n/c \rfloor$ for all integer $c>1$. 

\section{Algorithm} 

The starting point for the algorithm is a standard identity:

$$\sum_{x=1}^n M(\lfloor n/x \rfloor) = 1, n>1$$

Or, equivalently, 

$$M(n) = 1 - \sum_{x=2}^n M(\lfloor n/x \rfloor) $$ 

It's been noted in ~\cite{4} that the sequence $\lfloor n/x \rfloor$ takes on only $\lceil 2\sqrt{n} \rceil $ distinct values. We can rewrite the right-hand side as follows

$$M(n) = 1 - \sum_{x=2}^{n/t} M(\lfloor n/x \rfloor) - \sum_{x=1}^{t} (\lfloor \frac{n}{xt} \rfloor - \lfloor \frac{n}{(x+1)t} \rfloor )M(x)$$ 

%%%  check floors/ceils  %%%

where  $t > \lceil \sqrt{n} \rceil$. 

Now notice that the right hand side can be calculated recursively. In order to calculate M(n), we need to know M(n/2), M(n/3), M(n/4), .... In turn, M(n/2) requires M(n/4), M(n/6), etc. 

Suppose that we pick a value of $u>\sqrt{n}$ and create an array of $\lfloor n/u \rfloor$ elements, which would end up holding values of $M(\lfloor n/x \rfloor)$ for $1 \leq x \leq \lfloor n/u \rfloor$.

We then calculate all $M(y)$ for $1 \leq y\leq u$ in order by direct sieving (this would take $O(u^{1+\epsilon}$ time). For each $y$, we update the array to account for contributions of $M(y)$ according to the formula above. It's not hard to see that updating takes $O((n/u^{1/2})^{1+\epsilon})$ time. The overall time for these two tasks is minimized at $u=O(n^{2/3+\epsilon})$ and it also scales as $O(n^{2/3+\epsilon})$. 

Finally, when we are done sieving, it is a simple task to compute final values of M, and it takes negligible $O(n/u \ln{n/u}) \approx O(n^{1/3 + \epsilon})$ time.

Furthermore, it is easy to see that, if we want to compute $M(x)$ for multiple close values of $x$, we only need to do sieving once. Adjusting $u$ will result in optimal running time of $u=O(n^{2/3+\epsilon} N^{2/3+\epsilon})$, where $N$ is the number of values $x$ for which  we're trying to calculate $M(x)$.

The structure of the algorithm, with many large elements in the array that can be computed in parallel and more or less independently from each other, suggests that the algorithm may benefit from optimization for a general-purpose GPU (e.g. NVIDIA Fermi).

\section{Implementation} 

Since u is much larger than the amount of computer memory (for example, for $x=10^{22}$, $u \approx 10^{14}$, but modern desktop computers rarely have more than $\approx 10^{10}$ bytes of memory), it is necessary to divide calculations into blocks of 'y' and to handle each block in order.

For each block, the first step is to compute consecutive values of $M(y)$ throughout by sieving. Duplicate operations are avoided, and optimal performance is achieved if the block is greater than $u^{1/2}$. The nature of this computation is such that it involves a large number of random memory accesses across the entire block, and this is not a task that's particularly suitable for GPU optimization (among other reasons, the CPU has much larger L1 and L2 caches). Therefore, this task is done on the CPU.

A naive algorithm that performs sieving in a block $[y_1..y_2]$ can be described as follows:

1. Create an array A of integers with $y_2-y_1+1$ elements, corresponding to $y_1$, $y_1+1$, ..., $y_2$. Set these all to 1 initially.

2. For every prime number $p$ between 2 and $\lceil \sqrt{y_2} \rceil$, multiply elements in the array corresponding to multiples of $p$ by $-p$.

3. For every such prime, set all elements in the array corresponding to multiples of $p^2$ to 0.

4. Go through the array. For every element $y$, calculate the Mobius function $\mu(y)$:

 - If $A[y-y_1]$ is 0, then $\mu(y)$ is 0

 - Otherwise, if $|A[y-y_1]|$ is equal to $y$, then $\mu(y)$ is equal to the sign of $A[y-y_1]$

 - Otherwise, $\mu(y)$ is equal to minus the sign of $A[y-y_1]$

5. Add up all values of $\mu(y)$ to calculate $M(y)$.

The number of operations, according to the prime number theorem, is $O((y_2-y_1) \ln{\ln{\sqrt{y_2}}})$. While there's not much that can be done about the overall growth rate (obviously, the task can't be done faster than $O(y_2-y_1)$, and the log factor only plays a minor role, since, even  for $y_2 = 10^{20}$,  $\ln{\ln{\sqrt{y_2}}}) \approx 3.1$), we can improve the constant factor. One significant improvement can be achieved by modifying the temporary array to store approximate logarithms of prime products instead of exact values. This serves two purposes: it replaces each integer multiplication with an integer addition (which is often substantially faster); and it reduces the memory throughput of the algorithm by a factor of 8, because the array A in the range of values we're interested in has to contain 64-bit values, but base-2 logarithms can be stored in 8-bit variables. This is the modified algorithm:

0. (Preparation, has to be done once) Create an array of log-primes: For each prime $p_i$, the corresponding log-prime $l_i = \lfloor \log_2{p_i}\rfloor \operatorname{|} 1 $ (here '$|$' is binary "OR".)

1. Create an array A of 8-bit integers with $y_2-y_1+1$ elements, corresponding to $y_1$, $y_1+1$, ..., $y_2$. Set these all to 0 initially.

2. For every prime number $p$ between 2 and $\lceil \sqrt{y_2} \rceil$, add $l_i$ to all elements in the array corresponding to multiples of $p$.

3. For every such prime, set the most significant bit (by performing binary OR with 0x80) in all elements corresponding to multiples of $p^2$ within the block.

4. Go through the array. For every element $y$, calculate the Mobius function $\mu(y)$:

 - If the most significant bit of $A[y-y_1]$ is set, then $\mu(y)$ is 0
\
 - Otherwise, if $A[y-y_1]$ is greater than $\lfloor \log_2{y} \rfloor - 1$, then $\mu(y)$ is equal to $1-2(A[y-y_1] \operatorname{\&} 1)$ (here "\&" is binary "AND").

 - Otherwise, $\mu(y)$ is equal $-1+2(A[y-y_1] \operatorname{\&} 1)$

5. Add up all values of $\mu(y)$ to calculate $M(y)$.

It can be proved that sieving according to the algorithm will result in correct results, as long as $y_2$ is less than approximately $10^{27}$.

Further, we can precompute results of sieving with the first several primes and square primes, storing them in separate array, and copying this array instead of rerunning sieving. For example, the state of the temporary array after sieving with the first 4 primes and the first 2 square primes is periodic with the period of 2*2*3*3*5*7*11 = 13860.

Finally, the whole process can be modified to make use of multiple CPU cores, through the use of OpenMP. 

To perform the array updating in the minimal amount of time, we have to divide the (x,y) space into several distinct regions and process each region using a different procedure.

\begin{figure}
\centerline{\epsfig{file=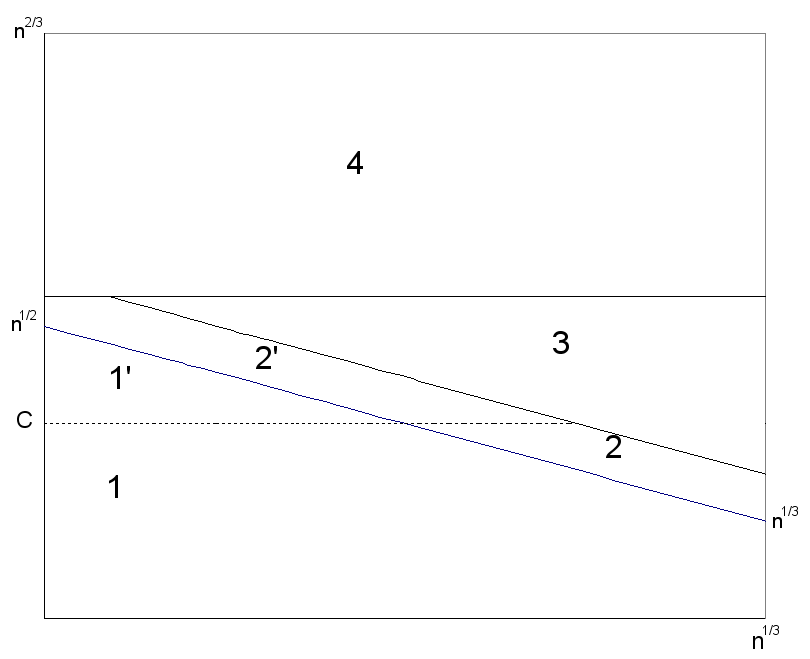,width=5in}}
\caption{Regions of (x,y) space}
\label{fig1}
\end{figure}

In regions 1 and 1' (defined as $y<c_1*(n/x)^{1/2}$ for some $c_1>1$), almost every (x,y) pair makes a nonzero contribution to one of the elements of the array. Adequate performance can be achieved by creating one GPU thread to handle a region with dimensions (1,512) and to add up contributions consecutively from every pair. 

An important factor that slows down computations in these regions is the necessity to perform at least one integer 64-bit division for every (x,y) pair. Integer division is an extremely slow operation that is not supported natively by the GPU and can take as much as a hundred clock cycles. (It's also quite slow on the CPU, but less so.) To improve performance, we can apply an algorithm such as ~\cite{5}, which accelerates integer division by replacing it with integer multiplication, addition, and shift. It requires that we precompute several constants for each value of the divisor, but, as long as each divisor is used multiple times, this results in net reduction of computation time. We can precompute all constants in advance for divisors up to some constant C, where C is only limited by the amount of available memory. This will result in substantial speedup as long as $C \gg u^{1/2}$. 

In regions 2 and 2'  (defined as $c_1*(n/x)^{1/2} \leq y<c_2*(n/x)^{1/2}$ where $c_2 \approx 10*c_1$),  there is a substantial number of (x,y) pairs which don't contribute anything, but the memory access pattern is relatively dense. Good performance is achieved when each GPU thread handles a region with dimensions (1,8192).

In the region 3, nonzero contributing pairs are sparse and processing time is limited by random memory access bandwidth. 

In the region 4 (defined as $y>c_3*n^{1/2}$, where $c \approx 2$), accesses are so sparse that it's not even justified to transfer results of the sieving onto the GPU any more: time spent transferring data over the PCI-e bus is larger than time that would be spent processing the same data with the CPU. At this point we no longer need the table of divisors used to speed up regions 1 and 2, so, if the amount of memory is a bottleneck, we can free that memory and increase the sieving block size up to several times $u^{1/2}$ (this is the only region where the ratio of block size to $u^{1/2}$ significantly greater than 1 is warranted).

\section{Approximate algorithm}

A different algorithm that can be employed in this task is based on the following remarkable formula:

$$M(x) \approx 2 \operatorname{Re} \sum_{i=1}^{\infty} \frac{x^{\rho_i}}{\rho_i \zeta'(\rho_i)} $$ (1)

or in a simple yet practical reformulation,

$$M(x)/x^{1/2} \approx q_n(x) = 2\sum_{i=1}^n a_i \cos{(z_i \ln{x} + b_i)}$$ (2)

where 

$$z_i = \operatorname{im} {\rho_i}$$
$$a_i = 1/|\rho_i \zeta'(\rho_i)|$$
$$b_i = -arg(\rho_i \zeta'(\rho_i))$$

and $\rho_i$ is the i'th zero of the Riemann zeta function. Somewhat more sophisticated reformulations, constructed by inserting a nontrivial kernel $f(n)$ inside the sum, are mentioned in the literature, but the formula above is particularly easy to calculate, and differences between different kernels for large $n$ are minor.

Although this formula is only approximate, it's usefulness here stems from the fact that we can get a good approximation - for example, two significant digits of $M(x)/x^{1/2}$ - with just 1000 terms, regardless of the magnitude of x. 

Calculating $q_n$ is also a task that nicely affords to GPU optimization. In this case, the implementation is relatively straightforward. The only nuance that needs addressing involves the use of single-precision vs. double-precision floating point. Modern GPUs are generally faster when operating in single precision than in double precision. This is particularly true with regard to the calculation of the cosine. All modern GPUs have native hardware support of basic single-precision transcendental functions such as the cosine. GPUs in NVIDIA Fermi family, for example, have the throughput of 4 or 8 single-precision cosines per clock cycle per multiprocessor: the NVIDIA GTX 560 can, in theory, exceed $10^{11}$ cosines per second (compare with modern desktop CPU's, which can only manage on the order of $10^8$ cosines/second). On the other hand, double-precision cosines are not natively supported by NVIDIA or AMD and they can be a factor of 20 or so slower. But the algorithm we have here can't be implemented completely in single precision, unless both n and x are small. 

\section{Results} 

To test the algorithm described in this article, a four-core Intel i5-2500 CPU running at 4.4 GHz with 16 GB of memory and a NVIDIA GeForce GTX 560 video card (384 processing cores at 1.7 GHz) with 1 GB of memory were used. Several choices of operating systems are available in this hardware configuration. Of these, x64 Linux (e.g. Ubuntu) is notable, because it comes with a C/C++ compiler (gcc/g++) that natively supports a 128-bit integer type. Since we're trying to compute M(n) for n above the 64-bit limit ($2^{64} \approx 1.9 * 10^{19}$), native 128-bit support substantially simplifies the code.

All the various constants and thresholds involved in this algorithm (such as dimensions handled by each GPU thread, or boundaries between regions 1/1' and 2/2') were tuned for peak performance. To verify correctness, many small values in $10^8 .. 10^{12}$ range were computed and results were compared with results obtained by naive O(n) time sieving. Values of Mertens function for $10^{13}, 10^{14}, 10^{15}$, and $10^{16}$ were compared against previously reported in the literature~\cite{4}. Calculating $M(10^{16})$ takes 1 minute: a substantial improvement, compared to the previously reported 1996 result where calculation for the same n took two days. (It should be noted that only a small part of this improvement is due to an overall increase in clock frequencies: the 1996 analysis employed a CPU running at 0.3 GHz, only one order of magnitude slower than hardware used here. Most of it is due to substantial parallelization and larger cache sizes.)

For higher n, direct verification is problematic, but results can be compared with approximate values. In order to make use of the approximate algorithm, the list of the first 2,000,000 zeros of Riemann zeta was downloaded from ~\cite{6}. Zeros were then refined to 30 digits and corresponding values of $\zeta'(\rho_i)$  were computed using Mathematica. To ensure correctness, these values for the first 200,000 zeros were independently computed to 30 digits using the arbitrary-precision floating point arithmetic library "mpmath" ver. 0.17 ~\cite{7}, and results were checked against each other to confirm the lack of major discrepancies. 

Approximations of M(n) for many large values of n were calculated. These turned out to be in agreement with exact calculations: the residual error between the exact and the approximate values of $M(n)/n^{1/2}$ for $n>~10^8$ was approximately normally distributed with the standard deviation of $\approx 4.2*10^{-4}$.

Finally, the code was used to compute values of M(n) for powers of 10 from $10^{17}$ to $10^{22}$. As expected, execution time scaled roughly as $n^{2/3}$. The computation for $10^{22}$ took 10 days. As described in the beginning of the article, computing $M(10^{22})$ simultaneously yields values of $M(y)$ for integer $y=\lfloor 10^{22}/c \rfloor$ for all integer c, including values such as $10^{19}$ and $10^{20}$; these agreed with their direct computations. Results are listed in table 1.

\begin{center}
Table 1. Values of $M(n)$ for select values of $n$.

\begin{tabular}[h]{|c|c|c|c|}\hline
$n$ & $M(n)$ & $M(n)/n^{1/2}$ & time \\\hline
$10^{16}$ & $-3195437$ & -0.032 & 1 m\\
$10^{17}$ & $-21830254$ & -0.069 & 4.5 m\\
$10^{18}$ & $-46758740$ & -0.047 & 20.7 m\\
$10^{19}$ & $899990187$ & 0.285 &  \\
$10^{20}$ & $461113106$ & 0.046 &  \\
$10^{21}$ & $3395895277$ & 0.107 &  \\
$10^{22}$ & $-2061910120$ & -0.021 & 10 d \\
\hline
\end{tabular}
\end{center}

    \begin{figure}
    \centerline{\epsfig{file=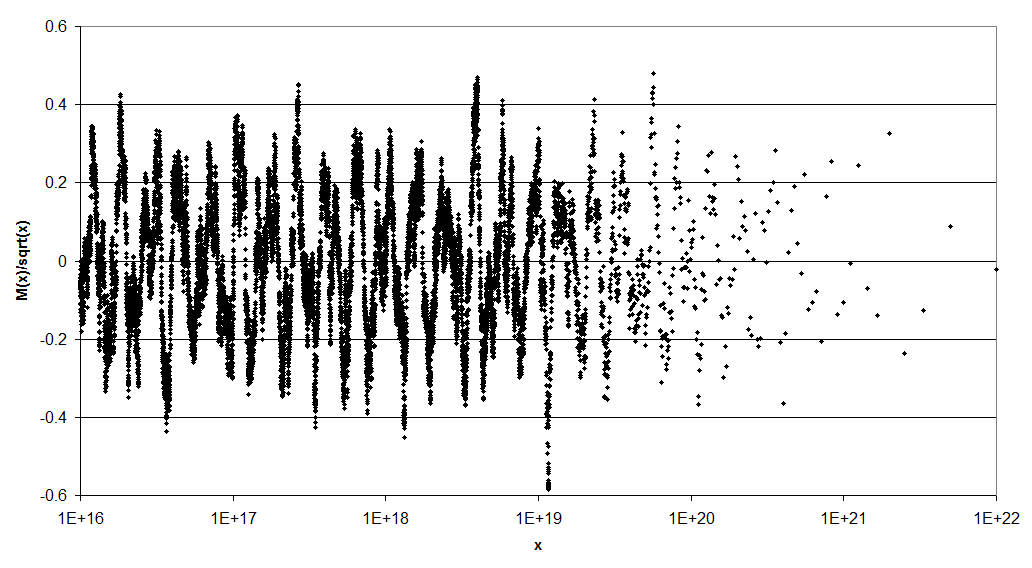,width=7in}} 
    \caption{ Values of M(x) for x between $10^{16}$ and $10^{22}$}
    \label{fig2}
    \end{figure} 

\section{New extreme value} 

As mentioned above, the largest value of $|M(n)/n^{1/2}|$ reported in the literature so far is 0.570591 at n=7766842813. ~\cite{3} analyzed approximate values for n up to $10^{10^{10}}$ using the Riemann zeta algorithm and noticed several larger potential extremes, including one potentially within reach at the present level of technology. Table 4 in that article notes that the $q_{10^6}$ approximation (based on the first $10^6$ zeros) reaches the value of -0.585 at $n \approx 1.16*10^{19}$, indicating a new extreme in the neighborhood.

To localize the extremum, a grid of values of M(n) in the neighborhood was systematically calculated. The density of the grid was increased till it became practical to cover intervals between the points of the grid with direct sieving (for example, if values of M(n) at $1.160986*10^{19}$ and $1.160987*10^{19}$ are known exactly, all values in the interval can be calculated, and the local extremum in this interval located, in a matter of hours). A plot with all computed points is shown in figure 2. 

The largest value of $|M(n)/n^{1/2}|$ found by this method turns out to be 0.585767684 (M(n)=-1,995,900,927) at $n=11,609,864,264,058,592,345$. 

Given the approach and the density of the grid, it is not possible to state with absolute certainty that the value is largest for, say, $n<10^{20}$, or even in the neighborhood of $1.16*10^{19}$. However, a heuristic argument can be made for negligible probability of finding a larger value below $10^{40}$.

Start by observing that a larger extreme value requires having a large value of $|p_{10^6}|$. Below $10^{40}$, the function $|p_{10^6}|$ only exceeds 0.575 inside the interval from $1.1600*10^{19}$ to $1.1612*10^{19}$ (the next such occurrence is in the vicinity of $2.26*10^{42}$). Since, as noted above, the residual error $M(n)/n^{1/2}-q_{10^6}(n)$  is approximately normally distributed with the standard deviation of $\approx 5.7*10^{-4}$, and the difference between 0.585767 and 0.575 comes out to 19 standard deviations, the probability for $|M(n)/n^{1/2}|$ at any randomly-picked value of n to exceed 0.585767 conditional on $q_{10^6}(n)$ being less than 0.575 is vanishingly small (less than $10^{-78}$) and therefore the probability for this to happen at any n below $10^{40}$ is still vanishingly small. (This, however, is subject to a caveat, which is mentioned in the next section.) So it remains to explore the interval of $1.1600*10^{19}$ to $1.1612*10^{19}$.

Suppose that we know values of $M(n)$ for $n=a$ and $n=b$. What can be said about the probability for M(n) to exceed a certan $M_0$ in the interval between $a$ and $b$, assuming that $b>a$, $M_0>M(a)$, and $M_0>M(b)$ (or, conversely, to go below $M_0$ if $M_0<M(a)$ and $M_0<M(b)$)? Locally, $M(n)$ can be thought of as essentially a random walk, which moves either up or down with the probability of $\sigma/2$ or sideways with the probability of $1-\sigma$ (where $\sigma=6/\pi^2$ is the share of squarefree numbers among all integers). The probability $P(M(n)>M_0)$ for such a random walk can be shown to be equal to 

$$P(M(n)>M_0) = \exp{(-(M(a)-M_0)(M(b)-M_0)/(2\sigma(b-a))}$$

to the first order in $b-a$ (which is sufficient for our needs). 

Of course, $M(n)$ is not a true random walk. It's not hard to see that a true random walk as described above would have values of $M(n)/n^{1/2}$ normally distributed with the standard deviation $6/\pi^2 \approx 0.608$. In reality, the standard deviation of $M(n)/n^{1/2}$ among computed values is only 0.11, suggesting that its dynamics involves a degree of regression towards the mean, and the probability estimate we have may be too high.

Nevertheless, using this this estimate and applying to all intervals shown in figure 2 yields that the total probability to find a larger extreme for $|M(n)/n^{1/2}|$ in the neighborhood (and, consequently, below $10^{40}$) is less than $0.05$.

\section{Search for the first counterexample} 

As the discussion above shows, it is highly improbable to find a counterexample to Mertens conjecture in the range where computation by the exact algorithm is tractable. Locating the first counterexample, therefore, is best approached using the approximate algorithm. 

Strictly speaking, the approximation (2) is only valid for all x in the assumption that the Riemann hypothesis is valid. If it fails, $|M(x)|$ can be expected to exceed $|x|$ for $x \approx \min_i t_i^{\frac{1}{s_i-0.5}}$, where $\rho_i = s_i + i t_i$ are nontrivial roots of Riemann zeta. Current limits on $\rho_i$ make it unlikely to see a counterexample below $10^{26}$ or so.

On the other hand, if the Riemann hypothesis holds, simple investigation of (2) leads to the conclusion that the conjecture is unlikely to be violated over a much wider range of values. A general sum of the form (2) can be treated as a sum of a large number of random variables. According to the central limit theorem, such a sum should be expected approximately normally distributed. A sum of the first $10^6$ terms ($q_{10^6}$) would be normally distributed with standard deviation of $\sqrt{2}\sum_{n=1}^{10^6} a_i^2 \approx 0.17$. In practice, the presence of several large non-normal terms in the sum means that the actual tail distribution of values of ($q_{10^6}$) is even narrower. 

The rest of this section assumes that the Riemann hypothesis holds.

The most straightforward algorithm that can be applied here is described in ~\cite{3}. It involves computing values of $p_n$ over a grid of values for some small value of $n$, recording values where $|q_n|$ exceeds a certain predefined threshold, and then studying areas where this happens in more detail. Unfortunately, this kind of brute-force search is too slow to yield the counterexample. Numerical evaluation of the probability distribution function (PDF) of (2) indicates that the search might need to be extended as high as $\exp{(10^{18})}$. Covering this region does not seem feasible at the present level of desktop PC technology: it would take an optimized software program and a typical high-end desktop PC hundreds of computer-years to reach $10^{18}$. It follows that a more intelligent algorithm needs to be used.

One useful observation that can be made is that the range of values that need to be checked can be constrained significantly just by examining the first several terms of (2). For example, the analysis of the PDF of (2) and its subsequences shows that any counterexample is unlikely to have $|q_7(x)|<0.425$, which is true for fewer than 1 in 10,000 values. 

In addition, all subsequences of (2) are quasiperiodic. This is not particularly useful for long subsequences, but for short sequences, such as $q_4$ and $q_7-q_4$, we can find $n<10^9$ such that these functions are "almost" periodic with periods $2\pi n/z_1$. For example,

$$(z_2/z_1)*274243136 = 407871396.0012$$
$$(z_3/z_1)*274243136 = 485262780.0010$$
$$(z_4/z_1)*274243136 = 590306027.0009$$
$$(z_5/z_1)*242101728 = 564116758.0001$$
$$(z_6/z_1)*242101728 = 643781791.9996$$
$$(z_7/z_1)*242101728 = 700862060.0012$$

which means that $q_4$ is approximately periodic with the period $274243136(2\pi/z_1)$ and $q_7-q_4$ is approximately periodic with the period $242101728(2\pi/z_1)$. We can accelerate the algorithm by computing values of $q_4$ and $q_7-q_4$ over one such period and then reusing them for multiple consecutive periods. Furthermore, we don't need to save actual values: it is much faster to compare each value against a predefined threshold and only save boolean flags which indicate whether the threshold is exceeded. This method ensures compact packing of the results (8 flags per byte) and minimizes memory troughput. 

Unfortunately, these algorithmic improvements are still insufficient to ensure locating a counterexample, but they allow us to extend coverage by two orders of magnitude.

To work with $x \gg \exp{(10^{10})}$, it is necessary to control for the loss of precision at high x. For example, computing $p_{10,000}(\exp{(10^{11})})$ involves computing a cosine of $z_{10,000}*10^{11} + b_{10,000} \approx 10^{15}$. Since these computations are done by the GPU, whose internal representation of double-precision floating point values has the machine precision $2^{-53} \approx 1.1*10^{-16}$, the computed value of the cosine could have be off by as much as 0.1. For much larger x or n, values of cosines could be essentially meaningless. To address this, a reference table of constants $a_i$, $b_i$ and $z_i$ to 30 significant digits was kept; for actual calculations, the range of study was divided into blocks of length $\approx 2*10^{11}$, and a new, 'shifted' table was generated for the every block using high-precision software according to the formula

$b'_i = (b_i+z_i*x_0) mod 2\pi$

The search was extended to $\exp{(10^{15})}$. Results up to $\exp{(10^{13})}$ were double-checked using brute-force search, and further verified using the 'mpmath' arbitrary-precision computing package. The search up to $\exp{(10^{15})}$ took approximately 7 days using the hardware mentioned in section 5. In line with expectations, no counterexamples were found. The largest value of $|q_{10^6}(x)|$ observed was $0.9565$ for $x \approx \exp{(5.0586*10^{14})}$.

Some "increasingly large positive/negative values" (using the terminology of ~\cite{3}) are listed in tables 2 and 3.

\begin{center}
Table 2. Increasingly large positive values.

\begin{tabular}[h]{|c|c|}\hline
$\ln{n}$ & $q_{10^6}(n)$ \\\hline
$22.773133$ & $+0.56959$\\
$97.526523$ & $+0.61863$\\
$984.282019$ & $+0.62512$\\
$1625.698493$ & $+0.62687$\\
$1957.803133$ & $+0.62849$\\
$2709.485814$ & $+0.65467$\\
$2794.384965$ & $+0.65955$\\
$12277.362671$ & $+0.79344$\\
$86458087.131684$ & $+0.80443$\\
$249548703.533702$ & $+0.81472$\\
$1467573228.501077$ & $+0.81702$\\
$1901564582.121964$ & $+0.82335$\\
$2500922487.505913$ & $+0.82884$\\
$3847517705.646364$ & $+0.83682$\\
$10407545552.85608$ & $+0.842485$\\
$21334043144.02927$ & $+0.86622$\\
$187114096628.77484$ & $+0.88636$\\
$1354137181464.62097$ & $+0.90578$\\
$6984497047106.74600$ & $+0.90744$\\
$84594507546024.46719$ & $+0.92208$\\
$117239588213313.90075$ & $+0.94102$\\
\hline
\end{tabular}
\end{center}
\begin{center}
Table 3. Increasingly large negative values.

\begin{tabular}[h]{|c|c|}\hline
$\ln{n}$ & $q_{10^6}(n)$ \\\hline
$43.898387$ & $-0.58478$\\
$140.373835$ & $-0.58940$\\
$853.851589$ & $-0.67715$\\
$3005.762748$ & $-0.68878$\\
$102494.024866$ & $-0.69580$\\
$150020.464414$ & $-0.70773$\\
$178259.151801$ & $-0.71541$\\
$203205.659988$ & $-0.74083$\\
$860440.495719$ & $-0.75254$\\
$1365643.292004$ & $-0.75927$\\
$2765781.628095$ & $-0.76041$\\
$7078384.260482$ & $-0.76879$\\
$13670267.747472$ & $-0.78505$\\
$19371574.223934$ & $-0.78747$\\
$57334128.09084$ & $-0.80765$\\
$167211796.14902$ & $-0.82488$\\
$405441986.398094$ & $-0.84497$\\
$4016980126.87193$ & $-0.85146$\\
$86339883457.03526$ & $-0.87296$\\
$264421251554.46918$ & $-0.89290$\\
$1278282683343.76520$ & $-0.89907$\\
$3680547202477.03623$ & $-0.91392$\\
$18747209824980.73961$ & $-0.91405$\\
$55714219637174.49540$ & $-0.92025$\\
$117892199597999.02070$ & $-0.92078$\\
$143697547951999.01914$ & $-0.93038$\\
$258592103887306.71643$ & $-0.94097$\\
$505863698785929.24318$ & $-0.95652$\\
\hline
\end{tabular}
\end{center}

\begin{figure}
\centerline{\epsfig{file=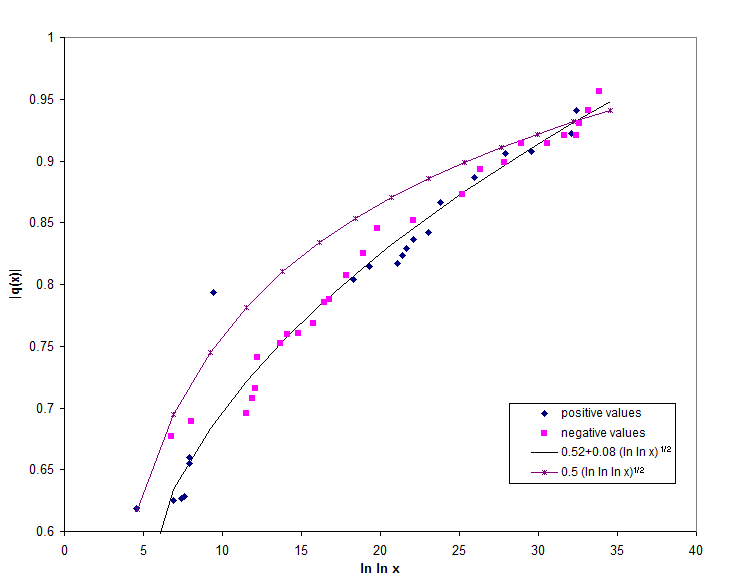}}
\caption{Increasingly large positive and negative values of $q(x)$}
\label{fig3}
\end{figure}

The distribution of values found here is in general agreement with expectations. Even though some authors (e.g. \cite{1}) observed tentative evidence for sign asymmetry (namely, that large positive values seemed to occur more often than large negative values), that does not appear to be the case here - the distribution found in this study is largely symmetric (see table 4.) 

Although no counterexamples were found in the area of study, further improvements to the algorithm and the computer hardware may bring this task within reach.

\begin{center}
Table 4. Tail ends of the frequency distribution of $|p_{7000}(x)|$  for $\ln{x} < 10^{15}$, sampled with the step $\pi/3072 z_1 \approx 7.24*10^{-5}$.

\begin{tabular}[h]{|c|c|c|}\hline
$|p_{7000}(x)|$ & positive values & negative values \\ \hline
$0.88 .. 0.89$ & $4499$  & $4399$ \\
$0.89 .. 0.90$ & $1472$  & $1518$ \\
$0.90 .. 0.91$ & $346$  & $444$ \\
$0.91 .. 0.92$ & $91$  & $112$ \\
$0.92 .. 0.93$ & $22$  & $42$ \\
$0.93 .. 0.94$ & $0$  & $6$ \\
$0.94 .. 0.95$ & $0$  & $6$\\
\hline
\end{tabular}
\end{center}
% 

%\medskip

%\noindent

\bigskip
%\newpage

 Eugene Kuznetsov (nameless@fastmail.fm)

\end{document}